\documentclass[10pt, conference, compsocconf]{IEEEtran}

\usepackage{amssymb,amsfonts,amsmath}

\usepackage[mathscr]{eucal}

\usepackage{stmaryrd}

\usepackage{amsbsy}

\usepackage{bm}

\usepackage{braket}

\newcommand{\Tra}{{\sf T}}

\newcommand{\qtext}[1]{\quad\text{#1}\quad}

\newcommand{\M}[2][]{{\bm{#1\mathbf{\MakeUppercase{#2}}}}} 
\newcommand{\Mn}[3][]{{\bm{#1\mathbf{\MakeUppercase{#2}}}}^{(#3)}} 
\newcommand{\MnTra}[3][]{{\bm{#1\mathbf{\MakeUppercase{#2}}}}^{(#3)\Tra}} 
\newcommand{\T}[2][]{\boldsymbol{#1\mathscr{\MakeUppercase{#2}}}} 
\newcommand{\Mz}[3][]{\M[#1]{#2}_{(#3)}}

\newcommand{\Tsize}[2]{\ensuremath#1_1 \times #1_2 \times \cdots \times #1_#2}

\newcommand{\Tentry}[2]{\ensuremath(#1_1,#1_2,\dots,#1_#2)}
\newcommand{\nrank}[1][n]{\ensuremath\text{rank}_#1}
\usepackage[caption=false,font=footnotesize,justification=centering]{subfig}

\usepackage{setspace}

\usepackage{xcolor}

\usepackage{algorithm,algpseudocode}

\makeatletter
\algrenewcommand\ALG@beginalgorithmic{\footnotesize}
\makeatother

\let\OLDthebibliography\thebibliography
\renewcommand\thebibliography[1]{
  \OLDthebibliography{#1}
  \setlength{\parskip}{0pt}
  \setlength{\itemsep}{0pt plus 0.3ex}
  \footnotesize
}

\usepackage{tikz}
\usetikzlibrary{3d}
\usetikzlibrary{patterns}
\usetikzlibrary{calc}
\usetikzlibrary{arrows}

\usepackage{pgfplots}
\usepackage{pgfplotstable}
\usepackage{etoolbox}
\usetikzlibrary{decorations.pathreplacing}
\makeatletter
\newcommand{\gettikzxy}[3]{%
  \tikz@scan@one@point\pgfutil@firstofone#1\relax
  \edef#2{\the\pgf@x}%
  \edef#3{\the\pgf@y}%
}
\makeatother

\usepackage{multirow}

\usepackage[naturalnames=false,hypertexnames=false]{hyperref}
\hypersetup{
  colorlinks = true,
  allcolors = blue
}

\usepackage[nameinlink]{cleveref}
\crefname{algorithm}{Alg.\@}{Algs.\@}
\crefname{figure}{Fig.\@}{Figs.\@}
\crefname{table}{Tab.\@}{Tabs.\@}
\crefname{section}{Sec.\@}{Secs.\@}

\newcommand{\lt}{\left}
\newcommand{\rt}{\right}

\newcommand{\datafile}{}
\newenvironment{inlinemath}{$}{$}

\begin{document}

\title{Parallel Tensor Compression for Large-Scale Scientific Data}

\author{%
  \IEEEauthorblockN{Woody Austin}
  \IEEEauthorblockA{University of Texas, Austin, TX, USA \\
    Email: austinwn@cs.utexas.edu}
  \and
  \IEEEauthorblockN{Grey Ballard and Tamara G. Kolda}
  \IEEEauthorblockA{Sandia National Labs, Livermore, CA, USA\\
    Email: gmballa@sandia.gov, tgkolda@sandia.gov}%
}

\maketitle

\begin{abstract}
As parallel computing trends towards the exascale, scientific data produced by high-fidelity simulations are growing increasingly massive. 
For instance, a simulation on a three-dimensional spatial grid with 512 points per dimension that tracks 64 variables per grid point for 128 time steps yields 8~TB of data, assuming double precision.
By viewing the data as a dense five-way tensor, we can compute a Tucker decomposition to find inherent low-dimensional multilinear structure, achieving compression ratios of up to 5000 on real-world data sets with negligible loss in accuracy. 
So that we can operate on such massive data, we present the first-ever distributed-memory parallel implementation for the Tucker decomposition, whose key computations correspond to parallel linear algebra operations, albeit with nonstandard data layouts.
Our approach specifies a data distribution for tensors that avoids any tensor data redistribution, either locally or in parallel.
We provide accompanying analysis of the computation and communication costs of the algorithms. 
To demonstrate the compression and accuracy of the method, we apply our approach to real-world data sets from combustion science simulations. 
We also provide detailed performance results, including parallel performance in both weak and strong scaling experiments.
\end{abstract}

\begin{IEEEkeywords}
  Tucker tensor decomposition; compression
\end{IEEEkeywords}

\IEEEpeerreviewmaketitle

\section{Introduction}
\label{sec:intro}

Today's high-performance parallel computers enable large-scale, high-fidelity simulations of natural phenomena across scientific domains.
As the speed and quality of simulations increase, the amount of data produced is growing at a rate that is creating bottlenecks in the scientific process.
A posteriori analysis of the data requires dedicated storage devices and parallel clusters even for simple computations.
A primary goal of this work is to provide a compression technique for large-scale simulation data which enables much more efficient data storage, transfer, and analysis, thereby facilitating bigger and better science.

Scientific simulation data is naturally multidimensional, tracking different variables in space and time; see \cref{fig:natural}. 
As a prototypical example, we focus on simulation data from combustion science research.
In this domain, simulated phenomena tend to be bursty, with important activity occurring in subsets of the spatial grid, small points in time, or involving a subset of the quantities of interest, like chemical species or fluid velocities.
Thus, the data typically have low-dimensional multilinear structure allowing for compression.
We consider compression based on the Tucker decomposition for higher-order tensors, which is analogous to principal component analysis (PCA) or the truncated singular value decomposition (T-SVD) of two-way data. 
We describe the Tucker decomposition in more detail in \cref{sec:tucker}.

\begin{figure}%
  \centering%
  \subfloat[Natural five-way multiway structure of scientific data.]{%
    \label{fig:natural}%
\begin{tikzpicture}
\newcommand{\drawtensor}
{
\coordinate (XFrontLowerLeft) at (0,0);
\draw (XFrontLowerLeft) rectangle ++ (\ix,\iy); %
\node[scale=0.6,above] at (0.5*\ix,0.5*\iy) {Spatial};
\node[scale=0.6,below] at (0.5*\ix,0.5*\iy) {Grid};
\begin{scope}[shift={(XFrontLowerLeft)},canvas is zx plane at y=\iy,rotate=90]
  \draw (0,0) rectangle ++ (\ix,\iz); %
\end{scope}
\begin{scope}[shift={(XFrontLowerLeft)},canvas is zy plane at x=\ix,rotate=90]
  \draw (0,0) rectangle ++ (\iy,\iz); %
\end{scope}
}
\newcommand{\drawtensorb}
{
\coordinate (XFrontLowerLeft) at (0,0);
\draw (XFrontLowerLeft) rectangle ++ (\ix,\iy); %
\begin{scope}[shift={(XFrontLowerLeft)},canvas is zx plane at y=\iy,rotate=90]
  \draw (0,0) rectangle ++ (\ix,\iz); %
\end{scope}
\begin{scope}[shift={(XFrontLowerLeft)},canvas is zy plane at x=\ix,rotate=90]
  \draw (0,0) rectangle ++ (\iy,\iz); %
\end{scope}
}

  \def\ix{1} %
  \def\iy{1} %
  \def\iz{0.6} %

  \foreach \y in {0,...,2}
  \foreach \x in {0,...,1}
  {
    \begin{scope}[shift={($(\x*\ix,\y*\iy)+(\x*0.5,\y*0.5)$)}]
      \drawtensor
    \end{scope}
  }

  \draw (-0.25,-0.25) rectangle +(3.25,4.75);
  \node[above,scale=0.75] at (1.5,4.5) {$\leftarrow$ Variables $\rightarrow$};
  \node[above,scale=0.75,rotate=90] at (-0.25,2) {$\leftarrow$ Time $\rightarrow$};

\end{tikzpicture}
%
%
%
%
  }
  \subfloat[Compression rates as fidelity varies for 550~GB simulation dataset.]{
    \label{fig:example_result}
\begin{tikzpicture}[scale=0.75]
\renewcommand{\datafile}{compression.dat}
\begin{loglogaxis}[
	width=2.25in,
	height=2.75in,
	ylabel={Compression Ratio}, 
	xlabel={Normalized RMS Error},
        ylabel near ticks,
]

	\addplot[nodes near coords*={$\pgfmathprintnumber[fixed,precision=0]\myvalue$},
        visualization depends on={\thisrow{sp_50} \as \myvalue},
        mark=*] 
        table[x=Error, y=sp_50, meta=sp_50] {\datafile};
\end{loglogaxis}
\end{tikzpicture}
%
%
%
%
  }
  \caption{}
\end{figure}

Our main contribution in this work is a distributed-memory parallel algorithm and implementation for computing the Tucker decomposition of general dense tensors.
Compression rates for a 550~GB scientific simulation dataset using our method are shown in \cref{fig:example_result}.
To the best of our knowledge, ours is the first distributed-memory implementation of a parallel algorithm for computing a Tucker decomposition.
Related work for other decompositions and algorithmic kernels are discussed in \cref{sec:related}.

The algorithm works for dense tensors of any order (i.e., number of dimensions) and size, given adequate memory, e.g.,  three times the size of the data.
The algorithm is efficient because it casts local computations in terms of {BLAS3} routines to exploit optimized, architecture-specific kernels; the data distributions and corresponding parallel computations are designed to reduce interprocessor communication.
We present the data distribution, parallel kernels, and overall algorithm in \cref{sec:par_dist,sec:kernels,sec:algorithms}, along with accompanying analysis of the computation costs and memory requirements.

Using real-world simulation data from combustion science, we demonstrate the effectiveness of Tucker for compression in \cref{sec:combustion}.
In particular, we show that these data sets have inherent low-dimensional multilinear structure that can be exploited for compression.
We show that the Tucker tensor decomposition can reduce the data by 50--75\% with normalized root mean squared (RMS) errors less than $10^{-6}$, and by 99.9\% and more with normalized RMS errors less than $10^{-2}$.
Such high compression rates allow terabytes of data to be reduced to gigabytes or megabytes, enabling easy data transfer and sharing. Additionally, we can reconstruct small subsets of the data upon request, enabling efficient analysis on even a single laptop.

The implementation is scalable. Results in \cref{sec:performance} show that the algorithm performs well for large and small data sets using up to $30{,}000$ cores.
It achieves near peak performance, as high as 66\%, on a single node consisting of 24 cores and up to 17\% of peak on over 1000 nodes.
At large scale, we are able to compress a 15~TB synthetic data set to 1.5~GB in about a minute and a 12~GB synthetic data set to 1~MB in under a second, with aggregate performance up to 100 TFLOPS.

\section{Tucker Tensor Decomposition}
\label{sec:tucker}

\subsection{Tensor Notation and Operations}

Let $\T{X}$ be a real-valued tensor of size $\Tsize{I}{N}$. 
We define
\begin{displaymath}
  I = \prod_{n=1}^N I_n \qtext{and} \hat I_{n} = I / I_n \text{ for } n \in \set{1,\dots,N}
\end{displaymath}
to be the total number of data elements and that number divided by the length of mode $n$,
respectively.
The \emph{mode-$n$ unfolding} rearranges the elements of $\T{X}$ to
form an $I_n \times \hat I_n$ matrix and is denoted by $\Mz{X}{n}$. Tensor
element $\Tentry{i}{N}$ maps to matrix element $(i_n,j)$ where
$
  j = 1 + 
  \sum_{k=1}^{n-1} (i_k-1) \hat I_{k} +
  \sum_{k={n+1}}^N (i_k-1) (\hat I_{k}/I_n).
$
The \emph{$n$-rank} of the tensor $\T{X}$ is the column rank of $\Mz{X}{n}$, denoted $\nrank(\T{X})$.  
The \emph{norm} of a tensor is the square root of the sum of the
squares of the entries, i.e.,
$ \|\T{X}\| = \|\Mz{X}{1}\|_F$.

The \emph{mode-$n$ product}
of the tensor $\T{X}$ with a real-valued
matrix $\M{V}$ of size $I_n \times J$ is denoted $\T{X} \times_n
\M{V}$ and is of size $I_1 \times \cdots \times I_{n-1} \times J
\times I_{n+1} \times \cdots \times I_N$. %
This can be expressed in terms of unfolded tensors, i.e.,
\begin{displaymath}
  \T{Y} = \T{X} \times_n \M{V} 
  \quad \Leftrightarrow \quad 
  \Mz{Y}{n} = \M{V}\Mz{X}{n}.
\end{displaymath}
This is also known as the tensor-times-matrix (TTM) product.
The order of multiplications is irrelevant, i.e., 
$  \T{X} \times_m \M{W} \times_n \M{V} =
 \T{X} \times_n \M{V} \times_m \M{W} \text{ for } m \neq n.$
If we are multiplying by a sequence of matrices, then
we may use the shorthand
$\T{Y} = \T{X} \times \set{\Mn{V}{n}}$,
to indicate that we should multiply $\T{X}$ by each matrix in the set
along the corresponding mode.
Unless otherwise indicated, we assume $\set{\Mn{V}{n}}$ is indexed from $n=1,\dots,N$.

\subsection{Tucker Algorithm and Key Kernels}

\begin{figure}[htbp]
  \centering
  %
\begin{tikzpicture}[scale=0.5,namenode/.style={scale=.75}]

\def\ix{3} %
\def\iy{3} %
\def\iz{2.5} %

\def\corescale{1.75}
\def\rx{\ix/\corescale}
\def\ry{\iy/\corescale}
\def\rz{\iz/\corescale}

\coordinate (XFrontLowerLeft) at (0,0);
\draw (XFrontLowerLeft) rectangle ++ (\ix,\iy); %
\begin{scope}[shift={(XFrontLowerLeft)},canvas is zx plane at y=\iy,rotate=90]
  \draw (0,0) rectangle ++ (\ix,\iz); %
\end{scope}
\begin{scope}[shift={(XFrontLowerLeft)},canvas is zy plane at x=\ix,rotate=90]
  \draw (0,0) rectangle ++ (\iy,\iz); %
\end{scope}
\node[namenode] at ($(XFrontLowerLeft) + (0.5*\ix, 0.5*\iy)$)  {$\T{X}$};

\coordinate (ApproxCtr) at ($(XFrontLowerLeft) + (\ix+0.4*\iz,0.75*\iy) + (0.75,0)$);
\node[namenode] at (ApproxCtr) {$\approx$};

\coordinate (U1LowerLeft) at ($(ApproxCtr) - (0,0.75*\iy) + (0.75,0)$);
\draw (U1LowerLeft) rectangle ++ (\ry,\iy);
\node[namenode] at ($(U1LowerLeft)+(0.5*\ry, 0.5*\iy)$)  {$\Mn{U}{1}$};

\coordinate (GFrontLowerLeft) at ($(U1LowerLeft) + (\ry+0.5,1)$);
\draw (GFrontLowerLeft) rectangle ++ (\rx,\ry);
\begin{scope}[shift={(GFrontLowerLeft)},canvas is zx plane at y=\ry,rotate=90]
  \draw (0,0) rectangle ++ (\rx,\rz);
\end{scope}
\begin{scope}[shift={(GFrontLowerLeft)},canvas is zy plane at x=\rx,rotate=90]
  \draw (0,0) rectangle ++ (\ry,\rz);
\end{scope}
\node[namenode] at ($(GFrontLowerLeft)+(0.5*\rx,.5*\ry)$)  {$\T{G}$};

\coordinate (U2LowerLeft) at ($(GFrontLowerLeft) + (\rx+\rz*0.4+0.5,0.5)$);
\draw (U2LowerLeft) rectangle ++ (\ix,\rx); %
\node[namenode] at ($(U2LowerLeft)+(0.5*\ix,0.5*\rx)$)  {$\Mn{U}{2}$};

\coordinate (U3LowerLeft) at ($(GFrontLowerLeft) + (0.5,\ry+.8)$);
\begin{scope}[shift={(U3LowerLeft)},canvas is zx plane at y=0,rotate=90]
  \draw (0,0) rectangle ++ (\rz,\iz); %
\end{scope}
\node[namenode] at ($(U3LowerLeft)+(1.2,0.5)$) {$\Mn{U}{3}$};

\end{tikzpicture}

%
%
%
%
  \caption{Tucker decomposition for $N=3$.}
  \label{fig:tucker}
\end{figure}
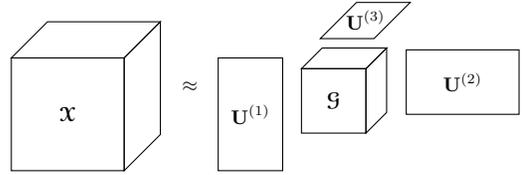

The Tucker decomposition \cite{Tu66} approximates a data tensor $\T{X}$ as
\begin{displaymath}
  \T{X} \approx \T{G} \times_1 \Mn{U}{1} \times_2 \Mn{U}{2} \cdots
  \times_N \Mn{U}{N} = \T{G} \times \set{ \Mn{U}{n} },
\end{displaymath}
where $\T{G}$ is the \emph{core tensor} of size $R_1 \times R_2 \times
\cdots \times R_N$ and $\Mn{U}{n}$ is a \emph{factor matrix} of size
$I_n \times R_n$ for $n=1,\dots,N$. 
We say $R_1 \times R_2 \times
\cdots \times R_N$ is the reduced dimension or, equivalently, the rank
of the reduced representation.
The decomposition is illustrated for $N=3$ in \cref{fig:tucker}.

Ideally, the data tensor has low-rank structure meaning that we
can choose $R_n \approx \nrank(\T{X}) \ll I_n$ so that most of the
variance of the data is preserved. 
Given factor matrices
$\set{\Mn{U}{n}}$, it is well known that
the optimal core is given by
$\T{G} = \T{X} \times \set{ \Mn{U}{n} }$ \cite{KoBa09}.

The storage is dominated by the core of size $R = \prod_{n=1}^N
R_n$. There is additional $\sum_{n=1}^N I_n R_n$ storage for the factor
matrices, but this is generally negligible compared to the storage of
the core.

The Tucker1 method, better known as truncated
higher-order singular value decomposition (T-HOSVD) \cite{DeDeVa00,Tu66}, is a
particular case of Tucker where $\Mn{U}{n}$ is set to be the $R_n$
leading left singular vectors of $\Mz{X}{n}$. 
The T-HOSVD is not optimal, but it is often a
good starting point for the iterative procedure described below.
We use a variation known as the \emph{sequentially-truncated HOSVD} (ST-HOSVD) \cite{VaVaMe12} for initialization. 
The first factor matrix is initialized as
for the T-HOSVD, i.e., the $R_1$ leading left singular values of
$\Mz{X}{1}$. The $n$th factor matrix is initialized as the $R_n$
leading left singular vectors of $\Mz{Y}{n}$ where 
$  \T{Y} = \T{X} \times_1 \MnTra{U}{1} \cdots \times_{n-1} \MnTra{U}{n-1}.$
The size of $\Mz{Y}{n}$ is $I_n \times ( \prod_{m<n} R_m)(\prod_{m>n} I_m)$.
One advantage of this method is that the $\T{Y}$ tensors are smaller
than $\T{X}$ for $n > 1$.
The ST-HOSVD is presented in \cref{alg:sthosvd}. 
Here, we pick the $R_n$ values according to a user-specified relative error threshold \cite{VaVaMe12}.
Although we do not make it explicit in the algorithm, the modes can be processed in any arbitrary order; see \cref{sec:vary_mode_order} for the impact of different orderings.

\begin{algorithm}[htb]
  \caption{Sequentially-Truncated HOSVD (ST-HOSVD)}
  \label{alg:sthosvd}
  \begin{algorithmic}[1]
    \Procedure{ST-HOSVD}{$\T{X}$, $\epsilon$}
    \State $\T{Y} \gets \T{X}$
    \For{$n=1,\dots,N$}\label{line:hosvd:loop}
    \State \label{line:hosvd:gram}
    $\M{S} \gets \Mz{Y}{n} \Mz{Y}{n}^{\Tra}$
    \State $R_n \gets $ min $R$ such that $\sum_{r>R} \lambda_r(\M{S}) \leq \epsilon^2\|\T{X}\|^2/N$
    \State \label{line:hosvd:evecs}
    $\Mn{U}{n} \gets $ leading $R_{n}$ eigenvectors of $\M{S}$
    \State \label{line:hosvd:ttm} $\T{Y} \gets \T{Y} \times_{n} \MnTra{U}{n}$
    \EndFor
    \State $\T{G} \gets \T{Y}$
    \State \Return $( \T{G},\set{\Mn{U}{n}} )$
    \EndProcedure
  \end{algorithmic}
\end{algorithm}

The higher-order orthogonal iteration (HOOI)
\cite{DeDeVa00a,KrDe80} is an alternating optimization method 
that further improves the approximation.  
The procedure cycles through the modes of
the tensor, calculating the leading left singular vectors of
$\Mz{Y}{n}$ where
$ \T{Y} = \T{X} \times \set{\MnTra{U}{m}}_{m \neq n},$
i.e., $\T{X}$ is multiplied in every mode except $n$ by the corresponding
factor matrix. The size of $\Mz{Y}{n}$ is
$I_n \times \hat R_n$ where $\hat R_n = \prod_{m \neq n} R_m$.
The HOOI method is presented in \cref{alg:hooi}.
HOOI is an iterative algorithm that monotonically improves the error but has no guarantees on converging to a global minimum; we iterate until the approximation error is small enough, the improvement in approximation falls below a given threshold, or a maximum number of iterations are reached.
We track the quantity $(\|\T{X}\|^2 - \|\T{G}\|^2)$  in
\cref{line:hooi:until} because it is equivalent to the fit of the
model, i.e., $\| \T{X} - \T{G} \times \set{\Mn{U}{n}} \|^2$ \cite{KoBa09}.

\begin{algorithm}[htb]
  \caption{Higher-order Orthogonal Iteration (HOOI)}
  \label{alg:hooi}
  \begin{algorithmic}[1]
    \Procedure{HOOI}{$\T{X}$, $\epsilon$}
    \State $(\T{G},\set{\Mn{U}{n}})=\text{ST-HOSVD}(\T{X},\epsilon)$ 
    \Repeat
    \For{$n=1,\dots,N$}
    \State \label{line:hooi:Y} 
    $\T{Y} \gets \T{X} \mathop{\times} \set{\MnTra{U}{m}}_{m \neq
      n}$
    \State \label{line:hooi:gram}
    $\M{S} \gets \Mz{Y}{n} \Mz{Y}{n}^{\Tra}$
    \State \label{line:hooi:evecs}
    $\Mn{U}{n} \gets $ leading $R_n$ eigenvectors of $\M{S}$
    \EndFor
    \State \label{line:hooi:core}
    $\T{G} \gets \T{Y} \times_N \MnTra{U}{N}$
    \Until{the quantity $(\|\T{X}\|^2 - \|\T{G}\|^2)$ ceases to
      decrease}\label{line:hooi:until}%
    \State \Return $(\T{G},\set{\Mn{U}{n}})$ 
    \EndProcedure
  \end{algorithmic}
\end{algorithm}

For both ST-HOSVD and HOOI, we compute the leading left singular
vectors of $\T{Y}$ by forming its $I_n \times I_n$ Gram matrix $\M{S} =
\Mz{Y}{n}\Mz{Y}{n}^{\Tra}$ and then computing its eigenvectors. 
Alternatively, we could work directly with $\Mz{Y}{n}$ and compute its singular vectors (see \cref{sec:conclusion} for more discussion). 
Our decision is motivated by the application-based assumption that $\epsilon$ is larger than the square root of machine precision and that $I_n$ is relatively small, i.e., $I_n \leq 2000$ for all $n$.

We focus on parallelizing the ST-HOSVD and HOOI methods whose inputs are $\T{X}$ and the desired accuracy, $\epsilon$. 
The three key operations in \cref{alg:sthosvd} and \cref{alg:hooi} are
\begin{enumerate}
\item the (sequence of) TTM operations to calculate $\T{Y}$,
\item the Gram matrix computations to calculate $\M{S}$, and
\item the eigenvector calculations to calculate $\Mn{U}{n}$.
\end{enumerate}
The difference between the two algorithms is that the tensors $\T{Y}$
are of different sizes. 
In HOOI, the core matrix is computed in
\cref{line:hooi:core}, exploiting the fact that the current $\T{Y}$
tensor already has the first $N-1$ products calculated. 

\subsection{Reconstruction}
\label{sec:reconstruction}

Given a core $\T{G}$ and factor matrices $\set{\Mn{U}{n}}$, we compute the approximate reconstruction of $\T{X}$, denoted $\T[\tilde]{X}$, by calling a sequence of TTM operations, i.e.,
\begin{equation}\label{eq:reconstruct}
  \T{X} \approx \T[\tilde]{X} = \T{G} \times \set{ \Mn{U}{n} }.
\end{equation}
Note that we can efficiently compute subtensors of $\T[\tilde]{X}$ (without forming the entire tensor) by modifying \cref{eq:reconstruct} appropriately, using subsets of rows of the factor matrices.

\section{Related Work}
\label{sec:related}

The Tucker decomposition is a powerful tool for compression of scientific data, as previously shown for hyperspectral images \cite{KaYaMe12} and volume rendering \cite{BaPa15}.
Our work parallelizes the method and demonstrates its utility on large-scale data sets.

To the best of our knowledge, ours is the first distributed memory implementation of the Tucker tensor decomposition. 
Zhou, Cichocki, and Xie~\cite{ZCX14} propose a randomized method for computing the Tucker decomposition of large tensors, but still assume all the data fits on a single machine.
Li et al.~\cite{LBPSV15} consider a shared-memory parallel implementation of dense TTM, which may be used directly and adapted for the local computations in our method. 

Several groups have parallelized the canonical polyadic (CP) tensor decomposition. 
The primary kernels are distinct from this work, but we mention a few similarities.
Karlsson, Kressner, and Uschmajew \cite{KKU14} have a parallel CP decomposition of multidimensional data with missing entries using either cyclic coordinate descent or alternating least squares; they use a similar data distribution strategy to ours in terms of the data tensor and factor matrices.
Kaya and U\c{c}ar \cite{KU15} focus on CP of \emph{sparse} tensors using a hypergraph partitioning scheme.
Kang et al.~\cite{KPHF12} provide a MapReduce implementation of CP for sparse data.
Smith et al.~\cite{SmRaSiKa15} parallelize a particular key kernel for CP on sparse data for three-mode tensors.
Phan and Cichocki \cite{PC11} break up the problem by computing CP decompositions of subtensors (in parallel) and then stitching the results together for a global CP decomposition. 

Another tensor decomposition known as the tensor train (TT) decomposition has recently been parallelized by Etter \cite{Etter15} by focusing on the recursive branching nature of the method.
There has also been considerable work on parallel tensor contraction; see, e.g., \cite{Martin} and references therein.

\section{Parallel Data Distributions}
\label{sec:par_dist}

For $N$-way tensors, we assume a logical $N$-way processor grid.
Let $\Tsize{P}{N}$ be the size of the processor grid.
For ease of presentation, we assume $P_n$ evenly divides $I_n$ and $R_n$, but our implementation does not require this.
Let $P = \prod P_n$  be the total number of processors and $\hat P_n = P / P_n$ be the number of processors in all modes but $n$.
From the point of view of a particular processor, we denote its local portion of a distributed object with an overhead bar.

\subsection{Tensor Distribution}
\label{sec:dist-tensor}
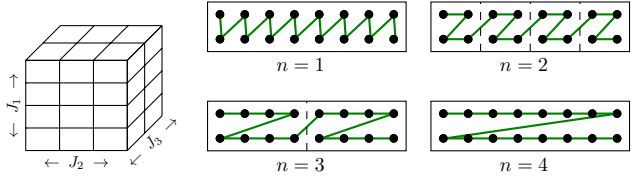
\begin{figure}[t]
\centering
\subfloat[Tensor  on $4 {\times} 3 {\times} 2$ processor grid.]{
\begin{tikzpicture}[scale=.3,textnode/.style={scale=.6}]

\def\ix{3} %
\def\iy{4} %
\def\iz{2} %
\def\sx{1.5} %
\def\sy{1}
\def\sz{2}
\def\offset{0.25}
\def\myquad{0.5}

\draw (0,0) grid[xscale=\sx,yscale=\sy] (\ix,\iy); %
\node[textnode,below] at (0.5*\ix*\sx,0) {$\leftarrow \; J_2 \; \rightarrow$};
\node[textnode,rotate=90,above] at (0, 0.5*\iy*\sy)  {$\leftarrow \; J_1 \; \rightarrow$};
\begin{scope}[shift={(XFrontLowerLeft)},canvas is zx plane at y=\iy*\sy,rotate=90]
  \draw (0,0) grid[xscale=\sx,yscale=\sz] (\ix,\iz); %
\end{scope}
\begin{scope}[shift={(XFrontLowerLeft)},canvas is zy plane at x=\ix*\sx,rotate=90]
  \draw (0,0) grid[xscale=\sy,yscale=\sz] (\iy,\iz); %
  \node[textnode,rotate=45,below] at (0,0.5*\iz*\sz) {$\leftarrow \; J_3 \; \rightarrow$};
\end{scope}
\end{tikzpicture}

%

%

%

%
%
%
%
%
%
%
%
%
%
%
%

%

%
%
%
%
%
%
%
%
%
%
%
%
%
%
%

%
%
%
%
%
\label{fig:tensor_block_dist}
}
\subfloat[Local data layout for $2 {\times} 2 {\times} 2 {\times} 2$ tensor unfolded in different modes.]{
\begin{tikzpicture}[scale=0.33] 
\newcommand{\datalayoutpic}
{
  \newcount\cnt
  \cnt=1
  \foreach \c in {1,...,\dimc}
  \foreach \a in {1,...,\dima}
  \foreach \b in {1,...,\dimb}
  {
    \node[draw,fill=black,minimum size=\circsize,inner sep=0pt,circle] (\the\cnt) at 
    ($(0,\dima) + (0,1) + (\b,0) - (0,\a) + (\dimb*\c,0) - (\dimb,0)$) {};
    \global\advance\cnt by 1
  }

  \pgfmathparse{\dima*\dimb*\dimc-1}
  \edef\total{\pgfmathresult}
  \foreach \i in {1,...,\total}
  {    
    \pgfmathparse{\i+1}
    \coordinate (A) at (\i);
    \coordinate (B) at ($(\pgfmathresult) -(\circsize,0)$);
    \draw[color=green!50!black,thick] (A)--(B);
  }

  \draw (0.5,0.5) rectangle +(8,2);

  \foreach \c in {1,...,\dimc}
  \foreach \a in {1,...,\dima}
  \foreach \b in {1,...,\dimb}
  \node[draw,fill=black,minimum size=\circsize,inner sep=0pt,circle] at 
  ($(0,\dima) + (0,1) + (\b,0) - (0,\a) + (\dimb*\c,0) - (\dimb,0)$) {};

}

  \begin{scope}[shift={(0,0)}]
  \def\dima{2}
  \def\dimb{1}
  \def\dimc{8}
  \def\circsize{1mm}
  \datalayoutpic   
  \node[below,scale=0.75] at (4.25,.5) {$n=1$};
  \end{scope}
  \begin{scope}[shift={(9,0)}] 
    \def\dima{2}
  \def\dimb{2}
  \def\dimc{4}
  \def\circsize{1mm}
  \datalayoutpic
  \pgfmathparse{\dimc-1}
  \foreach \x in {1,...,\pgfmathresult}
  \draw[dashed] ($(0.5,0.5)+(\x*\dimb,0)$) -- +(0,\dima);
  \node[below,scale=0.75] at (4.25,.5) {$n=2$};
  \end{scope}
  \begin{scope}[shift={(0,-4)}]
  \def\dima{2}
  \def\dimb{4}
  \def\dimc{2}
  \def\circsize{1mm}
  \datalayoutpic
  \pgfmathparse{\dimc-1}
  \foreach \x in {1,...,\pgfmathresult}
  \draw[dashed] ($(0.5,0.5)+(\x*\dimb,0)$) -- +(0,\dima);
  \node[below,scale=0.75] at (4.25,.5) {$n=3$};
  \end{scope}
  \begin{scope}[shift={(9,-4)}]
  \def\dima{2}
  \def\dimb{8}
  \def\dimc{1}
  \def\circsize{1mm}
  \datalayoutpic
  \node[below,scale=0.75] at (4.25,.5) {$n=4$};
  \end{scope}

\end{tikzpicture}
%

%
%
%
%
  \label{fig:unfolded}
}
\caption{Tensor data distribution/layout examples.}
\end{figure}

Let $\Tsize{J}{N}$ be the size of a generic tensor $\T{Y}$ where, for our purposes,
$J_n \in \set{I_n, R_n}$ for all $n$. 
Let $J = \prod J_n$ be the
total size and $\hat J_n = J / J_n$ be the size in all modes but $n$.
We impose a Cartesian parallel distribution of the tensor across
processors, which we refer to as a \emph{block distribution}. 
Each processor owns a distinct subtensor of size $J_1/P_1 \times
\cdots \times J_N/P_N$, with $J/P$ entries.
A similar tensor distribution is used in \cite{ZCX14}.
See \cref{fig:tensor_block_dist} for an illustration of a 3D 
tensor block distributed over a $4 \times 3
\times 2$ processor grid. 
The local tensor, $\T[\bar]{Y}$, is stored so that its mode-1 unfolding is in column-major order.

\subsection{Matrix Distribution}
\label{sec:matrix-distribution}

Factor matrices are also distributed across processors.
Given a mode $n$, let $\M{V}$ be a generic matrix of size $K \times J$ where $(K,J) = (R_n,I_n)$ for decomposition or $(K,J) =(I_n,R_n)$ for reconstruction. 
We treat our processor grid as two-dimensional of size $P_n \times \hat P_n$.
We divide $\M{V}$ into $P_n$ \emph{column blocks} so that
\begin{inlinemath}
  \M{V} =
  \begin{bmatrix}
    \M{V}_{1} & \M{V}_{2} & \cdots & \M{V}_{P_n}
  \end{bmatrix}.
\end{inlinemath}
We distribute the matrix \emph{redundantly} on every processor column, i.e., $\hat P_n$ times.
Since these matrices are relatively small, the redundant storage is negligible.
More precisely, if $R_n<\hat I_n/\hat P_n$, then the local storage of a factor matrix does not exceed the size of the local tensor.
Processor $\Tentry{p}{N}$ owns block column $\M{V}_{p_n}$ of size $K \times J/P_n$. 
We assume the local matrices are stored in row-major order.

\subsection{Unfolded Tensor Distribution}
\label{sec:dist-unfolded}

Unfolding a tensor is a purely logical process and involves no data redistribution.
Given a block distribution of a tensor, the unfolded tensor (a matrix) is also block distributed across a 2D processor grid.
In other words, a tensor unfolded in mode $n$ has dimension $J_n \times \hat J_n$ and is
distributed over a $P_n \times \hat P_n$ processor grid.
Note that if $P_n=1$, then the unfolded matrix has a 1D column distribution across $P$ processors.

The local portion of the unfolded tensor is equivalent to unfolding the local tensor. 
Again, this unfolding is logical; no local data distribution is required.
The local unfolded tensor, $\Mz[\bar]{Y}{n}$, is stored in $\prod_{m>n} (I_m/P_m)$ block columns, with each block column of size $(I_n/P_n) \times \prod_{m<n} (I_m/P_m)$ stored in row-major order.
The local data layout is illustrated in \cref{fig:unfolded}.
The dots show the elements of the unfolded tensor arranged as a matrix connected by a green line in the order that they are stored in memory. 
When $n=1$, $\Mz[\bar]{Y}{1}$ is in column-major order;
and when $n=4$, $\Mz[\bar]{Y}{4}$ is in row-major order.
For the interior modes, the data is a series of row-major subblocks.
For $n=2$, there are 4 subblocks of size $2 \times 2$. 
For $n=3$, there are 2 subblocks of size $2 \times 4$.
In local computations, each subblock can be processed separately using BLAS subroutines.

\section{Parallel Kernels and Analysis}
\label{sec:kernels}

\subsection{Parallel Cost Model and Collectives}
\label{sec:model}

To analyze our algorithms, we use the $\alpha$-$\beta$-$\gamma$ model of distributed-memory parallel computation, assuming 
the time to send a message of size $W$ words between any two processors is $\alpha+W\beta$, where $\alpha$ is the latency cost, $\beta$ is the per-word transfer cost, and $\gamma$ is the time for one floating point operation (flop).
For more discussion of the model and descriptions of efficient collectives, see \cite{CH+07,TRG05}.
Let $W$ be the data size and $P$ be the number of processors, then the costs for relevant operations are summarized in \cref{tab:communication_model}.
For simplicity of presentation, we will ignore the flop cost of reduce and all-reduce, as they are typically dominated by the bandwidth costs.
\begin{table}[ht]
  \centering\footnotesize
  \caption{Communication costs in $\alpha$-$\beta$-$\gamma$ model.}
  \label{tab:communication_model}
  \begin{tabular}{|l|l|} \hline
    Send/Receive & $\alpha + \beta W$ \\
    All-gather & $\alpha \log P + \beta \frac{P-1}{P} W$ \\
    Reduce & $\alpha \log P + (\beta + \gamma) \frac{P-1}{P} W$ \\
    All-reduce & $2 \alpha \log P + (2\beta + \gamma) \frac{P-1}{P} W$ \\
    \hline
  \end{tabular}
\end{table}

\subsection{Parallel TTM Computation}
\label{sec:ttm}

\begin{figure*}[t]
\centering
\begin{tikzpicture}[%
  scale=0.75,%
  namenode/.style={scale=1},%
  sizenode/.style={scale=0.7,midway},%
  curlybrace/.style={decorate,decoration={brace,amplitude=2pt}},%
]

\colorlet{color1}{red!30}
\colorlet{color2}{blue!30}
\colorlet{color3}{green!30}

\def\nameyoffset{0.75}
\def\curlybraceoffset{0.1}
\def\arrowoffset{0.1}

\def\dima{1.5} %
\def\dimb{2} %
\def\dimc{8} %

\def\divb{3} %
\def\divc{4} %
\edef\diva{\divb} %
\def\sdivc{3} %

\pgfmathparse{\dimb/\divb}
\edef\incb{\pgfmathresult} %
\pgfmathparse{\dima/\diva}
\edef\inca{\pgfmathresult} %
\pgfmathparse{\dimc/\divc}
\edef\incc{\pgfmathresult} %

\coordinate (VLL) at (0.5,1.5); %
\begin{scope}[shift=(VLL)]
  \draw [xscale=\incb,yscale=\dima] grid(\divb,1);

  \foreach \x in {1,...,\divb}
    \draw [fill=color\x] ($(\x*\incb,0)-(\incb,0)$) rectangle +(\incb,\dima);

  \pgfmathparse{\diva-1}
  \foreach \x in {1,...,\pgfmathresult}
    \draw[dashed,thick] (0,\x*\inca) -- +(\dimb,0);

  \draw [curlybrace] (-\curlybraceoffset,0) -- +(0,\dima) node [sizenode,left] {$K$};

  \draw [curlybrace] (0,\dima+\curlybraceoffset) -- +(\incb,0) node [sizenode,above] {$J_n/P_n$};
\end{scope}

\gettikzxy{(VLL)}{\ptx}{\pty}
\node[namenode] at ($(\ptx,0)+(\dimb / 2,\nameyoffset)$) {$\M{V}$};

\coordinate (BULLET) at ($(VLL) + (\dimb, \dima / 2) + (0.5,0)$);
\node[namenode] at  (BULLET) {$\bullet$};

\coordinate (YLL) at ($(BULLET) + (0.5,0) - (0,\dimb/2)$); %
\begin{scope}[shift=(YLL)]
  \draw [xscale=\incc,yscale=\incb] grid(\divc,\divb);
  \foreach \x in {1,...,\divb}
    \draw [fill=color\x] ($(0,\dimb)-(0,\x*\incb)$) rectangle +(\incc,\incb);
  \draw [curlybrace] (0,\dimb+\curlybraceoffset) -- +(\incc,0) node [sizenode,above] {$\hat J_n/\hat P_n$};

  \draw [curlybrace] ($(-\curlybraceoffset,\dimb) - (0,\incb)$) -- +(0,\incb) node [sizenode,left] {$\displaystyle\frac{J_n}{P_n}$};

\end{scope}

\gettikzxy{(YLL)}{\ptx}{\pty}
\node[namenode] at ($(\ptx,0) + (\dimc/2,\nameyoffset)$) {$\Mz{Y}{n}$};

\coordinate (EQ) at ($(YLL)+(\dimc,\dimb/2)+(0.5,0)$);
\node[namenode] at (EQ) {$=$};

\coordinate (ZLL) at ($(EQ) + (.5,-\dima/2)$);

\gettikzxy{(ZLL)}{\ptx}{\pty}
\node[namenode] at ($(\ptx,0) + (\dimc/2,\nameyoffset)$) {$\Mz{Z}{n}$};

\begin{scope}[shift=(ZLL)]
  \draw [xscale=\incc,yscale=\inca] grid(\divc,\diva);

  \foreach \x in {1,...,\divb}
    \draw [fill=color\x] ($(0,\dima)-(0,\x*\inca)$) rectangle +(\incc,\inca);

  \draw [curlybrace] (0,\dima+\curlybraceoffset) -- +(\incc,0) node [sizenode,above] {$\hat J_n/\hat P_n$};

  \draw [curlybrace,decoration={mirror}] ($(\dimc,\dima) + (\curlybraceoffset,-\inca)$) -- +(0,\inca) node [sizenode,right] {$\displaystyle\frac{K}{P_n}$};

  \foreach \x in {1,...,\diva}
  \node[circle,scale=0.7] at ($(\incc/2,\dima)-(0,\x*\inca)+(0,\inca/2)$) {\x};
\end{scope}

\end{tikzpicture}

%
%
%
%
\caption{Parallel distributions of matrices involved in TTM
  operation $\T{Z} = \T{Y} \times_2 \M{V}$ on a $2 \times 3 \times 2$ processor grid. The data owned by 3 of the 12 processors is color coded. The computation requires $P_2=3$ iterations, each using one block row of $\M{V}$ and producing one (labeled) block of the output.}
\label{fig:par_ttm}
\end{figure*}
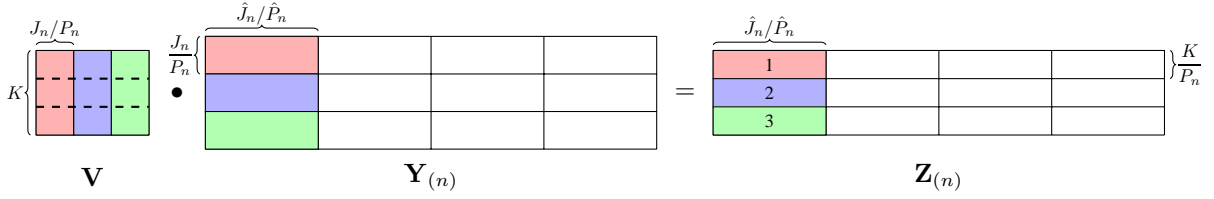

We consider the parallel algorithm for the TTM
operation,
\begin{displaymath}
  \T{Z} = \T{Y} \times_n \M{V},
\end{displaymath}
for a generic tensor $\T{Y}$ of size $\Tsize{J}{N}$ and a matrix
$\M{V}$ of size $J_n \times K$. The result is a tensor $\T{Z}$
of size $J_1 \times \cdots \times J_{n-1} \times K \times J_{n+1}
\times \cdots \times J_N$.
We can rewrite the operation in matricized form as
\begin{equation}
  \nonumber
  \Mz{Z}{n}=\M{V}\Mz{Y}{n}.
\end{equation}

The matrices are distributed as described in \cref{sec:par_dist}, and our parallel algorithm is presented in \cref{alg:TTM}.
We divide $\M[\bar]{V}$, the local portion of $\M{V}$, of size $K \times J_n/P_n$, 
into $P_n$ block rows, i.e.,
$\M[\bar]{V}^{[\ell]}$ denotes the $\ell$th block row of size $K/P_n \times J_n/P_n$.
We compute the local product one block row at a time, i.e., locally compute
\begin{displaymath}
\T{W} = \bar{\T{Y}} \times_n \M[\bar]{V}^{[\ell]} \quad \text{ or } \quad  \Mz{W}{n} = \M[\bar]{V}^{[\ell]} \Mz[\bar]{Y}{n},
\end{displaymath}
and sum the results across all processors in the same column to yield the result for the $\ell$th member.
The local computation is also a mode-$n$ TTM, which
we implement using \texttt{dgemm} within BLAS while respecting the local layout of the unfolded tensor as described in \cref{sec:dist-unfolded}.

The blocking ensures that the size of the intermediate products, which is $K/P_n \times \hat J_n / \hat P_n$, is never larger than the size of the local result tensor.
We note that if $K < J_n/P_n$, then we can avoid the blocking strategy, performing a single local matrix multiplication followed by a single reduce-scatter, without the temporary memory exceeding the size of the local input tensor.
This optimization reduces latency cost but does not affect computation or bandwidth costs; our analysis will assume the blocking strategy but our implementation exploits the non-blocked approach when possible.

\begin{algorithm}
  \caption{Parallel TTM}
  \label{alg:TTM}
  \begin{algorithmic}[1]
    \Procedure{TTM}{$\T{Y}$,$\M{V}$,$n$}
    \State $\texttt{myProcID} \gets \Tentry{p}{N}$
    \State $\texttt{myProcCol} \gets (p_1,\dots,p_{n-1},\ast,p_{n+1},\dots,p_N)$
    \For{$\ell=1,\dots,P_n$}
    \State $\T{W} \gets \bar{\T{Y}} \times_n \M[\bar]{V}^{[\ell]}$
    \label{line:TTM:matvec}
    \State $\T[\bar]{Z} \gets \textsc{Reduce}(\T{W},\texttt{myProcCol}, \ell)$
    \Comment Root is $p_n{=}\ell$
    \label{line:TTM:reduce}
    \EndFor
    \State \Return $\T{Z}$
    \EndProcedure
  \end{algorithmic}
\end{algorithm}

\Cref{fig:par_ttm} illustrates the computation on a $2 \times 3 \times
2$ processor grid for $n=2$. 
We implicitly treat the processor grid as $P_n \times \hat P_n = 3 \times 4$ and consider $\Mz{Y}{n}$ to be block distributed as described in \cref{sec:dist-unfolded}.  
We color code the block column of $\Mz{Y}{n}$ owned by the first column in the processor grid; these three processors work together to compute the first block column of the result.
Note that $\M{V}$ is redundantly stored across every processor column as described in \cref{sec:matrix-distribution}.
The result tensor $\T{Z}$ is partitioned in the same way as the input tensor $\T{Y}$.
Note that the blocks of the result are calculated one at a time, and numbered in the figure in the order that they are computed.

The cost and memory of the parallel TTM (\cref{alg:TTM}) is
\begin{align*}
  C_{\textsc{ttm}} =& 
  \underbrace{ 2\gamma \frac{JK}{P} }_{ P_n \times \text{ \cref{line:TTM:matvec}}} +
  \underbrace{ \alpha P_n \log P_n + \beta (P_n-1) \frac{\hat J_n K}{P} }_{P_n \times \text{ \cref{line:TTM:reduce}}}, \text{ and}\\
  M_{\textsc{ttm}} = &\underbrace{J/P}_{\T[\bar]{Y}} + 
  \underbrace{J_nK/P_n}_{\M[\bar]{V}} + 
  \underbrace{\hat J_n K/P}_{\T[\bar]{Z}} + 
  \underbrace{\hat J_n K / P}_{\T{W}} .
\end{align*}
The storage for $\T{W}$ is temporary.
If $P_n = 1$, then no parallel communication is required.

\subsection{Parallel Gram Computation}
\label{sec:gram}

Given mode $n$, we compute the Gram matrix $\M{S} = \Mz{Y}{n}\Mz{Y}{n}^\Tra$ where $\T{Y}$ is a tensor of size $J_1 \times \cdots \times J_N$. 
In the context of ST-HOSVD and HOOT, we know $J_n = I_n$.
The unfolded tensor  $\Mz{Y}{n}$ is block distributed on the $P_n \times \hat P_n$ processor grid per \cref{sec:dist-unfolded}.
The result $\M{S}$, of size $J_n \times J_n$, will be distributed with respect to mode $n$ as described in \cref{sec:matrix-distribution}.
We ignore the fact that $\M{S}$ is symmetric, storing both upper and lower triangles explicitly.

\begin{figure}
  \centering
\begin{tikzpicture}[%
  scale=0.70,%
  namenode/.style={scale=1},%
  sizenode/.style={scale=0.7,midway},%
  curlybrace/.style={decorate,decoration={brace,amplitude=2pt}},%
]
\colorlet{color1}{red!30}
\colorlet{color2}{blue!30}
\colorlet{color3}{green!30}

\def\nameyoffset{0.5}
\def\curlybraceoffset{0.1}
\def\arrowoffset{0.1}

\def\dima{2} %
\def\dimb{3.5} %
\def\diva{3} %
\def\divb{4} %

\pgfmathparse{\dima/\diva}
\edef\inca{\pgfmathresult} %
\pgfmathparse{\dimb/\divb}
\edef\incb{\pgfmathresult} %

\coordinate(YLL1) at ($(1,1)+(0,\dimb/2)-(0,\dima/2)$);
\begin{scope} [shift=(YLL1)]
  \draw [xscale=\incb,yscale=\inca] grid(\divb,\diva);

  \foreach \x in {1,...,\diva}
    \draw [fill=color\x] ($(0,\dima)-(0,\x*\inca)$) rectangle +(\incb,\inca);

  \draw [curlybrace] (0,\dima+\curlybraceoffset) -- +(\incb,0) node [sizenode,above] {$\hat J_n/\hat P_n$};

  \draw [curlybrace] ($(-\curlybraceoffset,\dima) - (0,\inca)$) -- +(0,\inca) node [sizenode,left] {$\displaystyle\frac{J_n}{P_n}$};
\end{scope}

\gettikzxy{(YLL1)}{\ptx}{\pty}
\node[namenode] at ($(\ptx,0) + (\dimb/2,\nameyoffset)$) {$\Mz{Y}{n}$};

\coordinate (BULLET) at ($(YLL1) + (\dimb, \dima / 2) + (0.5,0)$);
\node[namenode] at  (BULLET) {$\bullet$};

\coordinate(YLL2) at ($(BULLET) + (0.5,0) - (0,\dimb/2)$);
\begin{scope} [shift=(YLL2)]
  \draw [xscale=\inca,yscale=\incb] grid(\diva,\divb);

  \foreach \x in {1,...,\diva}
    \draw [fill=color\x] ($(-\inca,3*\incb)+(\x*\inca,0)$) rectangle +(\inca,\incb);

\end{scope}

\gettikzxy{(YLL2)}{\ptx}{\pty}
\node[namenode] at ($(\ptx,0) + (\dima/2,\nameyoffset)$) {$\Mz{Y}{n}^{\Tra}$};

\coordinate (EQ) at ($(YLL2)+(\dima,\dimb/2)+(0.5,0)$);
\node[namenode] at (EQ) {$=$};

\coordinate (SLL) at ($(EQ) + (.5,-\dima/2)$);

\def\Ssp{.25}
\coordinate (SLL1) at ($(SLL)+(\Ssp,\Ssp)$);
\coordinate (SLL2) at ($(SLL1)+(\Ssp,\Ssp)$);
\coordinate (SLL3) at ($(SLL2)+(\Ssp,\Ssp)$);
\draw [shift=(SLL3),xscale=\inca,yscale=\inca] grid(\diva,\diva);
\fill [white] (SLL2) rectangle ($(SLL2)+(\dima,\dima)$);
\draw [shift=(SLL2),xscale=\inca,yscale=\inca] grid(\diva,\diva);
\fill [white] (SLL1) rectangle ($(SLL1)+(\dima,\dima)$);
\draw [shift=(SLL1),xscale=\inca,yscale=\inca] grid(\diva,\diva);

\begin{scope} [shift=(SLL)]

  \foreach \x in {1,...,\diva}
    \draw [fill=color\x] ($(-\inca,0)+(\x*\inca,0)$) rectangle +(\inca,3*\inca);

  \draw [xscale=\inca,yscale=\inca] grid(\diva,\diva);

  \foreach \x in {1,...,\diva}
    \foreach \y in {1,...,\diva}
      \node[namenode] at ($(\x*\inca,\y*\inca)-(\inca/2,\inca/2)$) {+};

  \draw [curlybrace,decoration=mirror] (0,-\curlybraceoffset) -- +(\inca,0) node [sizenode,below] {$J_n/P_n$};

  \draw [curlybrace] ($(-\curlybraceoffset,\dima) - (0,3*\inca)$) -- +(0,\inca) node [sizenode,left] {$\displaystyle\frac{J_n}{P_n}$};
\end{scope}

\gettikzxy{(SLL)}{\ptx}{\pty}
\node[namenode] at ($(\ptx,0) + (\dima/2,\nameyoffset)$) {$\M{S}$};

\end{tikzpicture}
%
%
%
%
  \caption{Parallel distribution of matrices involved in Gram operation $S = \Mz{Y}{2} \Mz{Y}{2}^T$ on a $2 \times 3 \times 2$ processor grid. The data owned by 3 of the 24 processors is color coded. Each processor column computes local matrix-vector products and the results are summed across each block row.}
  \label{fig:gram}
\end{figure}
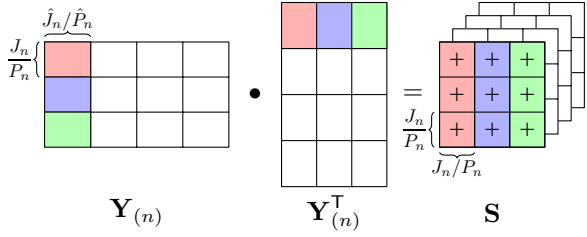

The data distribution for Gram is illustrated in \cref{fig:gram}, and the method is shown in \cref{alg:Gram}.
Each processor column owns a block column of $\Mz{Y}{n}$ and computes an intermediate matrix $\M{V}$ of dimension $J_n \times J_n/P_n$.
The matrix $\M{V}$ is computed in row blocks of dimension $J_n/P_n \times J_n/P_n$, where each row block is the product of an unfolding of the local tensor with the unfolding of another processor's local tensor.
The computation at \cref{line:Gram:multiply0} is a local Gram computation that we perform using \texttt{dsyrk} within {BLAS}, though for interior modes it requires multiple subroutine calls to respect the local layout (as in the case of TTM).
The computation at \cref{line:Gram:multiply} is a nonsymmetric analogue that we perform using \texttt{dgemm}.
The $\M{V}$ matrices are then summed across each processor row ($\hat P_n$ processors), using an all-reduce collective, so that the result (a block column of $\M{S}$) is replicated across the processor row.
Note that if $P_n=1$, then the computation fully exploits symmetry, and the only communication is the all-reduce across $P$ processors.

\begin{algorithm}
  \caption{Parallel Gram}
  \label{alg:Gram}
  \begin{algorithmic}[1]
    \Procedure{gram}{$\T{Y}$,$n$}
    \State $\texttt{myProcID} \gets \Tentry{p}{N}$
    \State $\texttt{myProcCol} \gets (p_1, \dots, p_{n-1}, \ast, p_{n+1}, \dots, p_N)$
    \State $\texttt{myProcRow} \gets (\ast, \dots, \ast, p_k, \ast, \dots, \ast)$
    \State $\M{V}^{[p_n]} \gets \Mz[\bar]{Y}{n} \Mz[\bar]{Y}{n}^\Tra$
    \label{line:Gram:multiply0}
    \For{$i=1$ to $P_n-1$}
    \State $j \gets (p_n-i) \mod P_n$ 
    \State $k \gets (p_n+i) \mod P_n$ 
    \State Send $\bar{\T{Y}}$ to process $(p_1, \dots, p_{n-1}, j, \dots, p_N)$ \label{line:Gram:send}
    \State  Receive $\T{W}$ from process $(p_1, \dots, p_{n-1}, k, \dots, p_N)$
\label{line:Gram:recv}
    \State $\M{V}^{[k]} \gets \Mz[\bar]{Y}{n} \Mz{W}{n}^{\Tra}$
    \label{line:Gram:multiply}
    \EndFor
    \label{line:Gram:end}
    \State $\M[\bar]{S} = \text{All-Reduce}(\M{V},\texttt{myProcRow})$
    \label{line:Gram:all-reduce}
    \State \Return $\M[\bar]{S}$
    \EndProcedure
  \end{algorithmic}
\end{algorithm}

The cost of parallel Gram (\cref{alg:Gram}) is
\begin{align*}
  C_{\textsc{gram}} = &
  \underbrace{ \gamma 2 {J_n J}/{P} }
  _{\text{\cref{line:Gram:multiply0} } + (P_n-1) \times \text{ \cref{line:Gram:multiply}}} 
  + \underbrace{2(P_n-1)\lt( \alpha + \beta {J}/{P} \rt) }
  _{(P_n-1) \times \text{\cref{line:Gram:send,line:Gram:recv}}} \\
  & +  \underbrace{2\alpha \log \hat P_n + 2\beta (\hat P_n-1)J_n^2 / P}
  _{\text{\cref{line:Gram:all-reduce}}}, \text{ and}\\
  M_{\textsc{gram}} = &
  \underbrace{J/P}_{\T[\bar]{Y}} + 
  \underbrace{J/P}_{\T{W}} + 
  \underbrace{J_n^2/P_n}_{\M{V}} + 
  \underbrace{J_n^2/P_n}_{\M[\bar]{S}}.  
\end{align*}
The storage for $\T{W}$ and $\M{V}$ is temporary.
Note that  up to a factor of two could be saved by exploiting symmetry of $\M{S}$.

\subsection{Parallel Eigenvectors Computation}
\label{sec:eigenvectors}

\Cref{alg:eigenvectors} presents our method for computing the leading eigenvectors of the Gram matrix.
After the Gram computation, the matrix $\M{S}$  of size $I_n \times I_n$ 
is stored redundantly on every processor column in the $P_n \times \hat P_n$ 
processor grid as described in \cref{sec:matrix-distribution}.
We enforce the same distribution of the transpose of the output $I_n \times R_n$ eigenvector matrix $\Mn{U}{n}$, which implies a block row distribution of $\Mn{U}{n}$.
Because we assume $I_n$ is relatively small, e.g., $I_n \leq 2000$, our approach is essentially a sequential algorithm. 
We perform an all-gather across the processor fiber so that every processor owns the entire matrix $\M{S}$. Every processor performs the local eigenvector computation redundantly using \texttt{dsyevx} within LAPACK and then extracts the appropriate subset of its local result to obtain the desired final distribution.
The algorithm is presented in \cref{alg:eigenvectors}.
While we assume $R_n$ is an input to the algorithm, we can also choose $R_n$ ``on the fly'' based on the desired error threshold for similar cost.

\begin{algorithm}
\caption{Parallel Eigenvectors}
\label{alg:eigenvectors}
\begin{algorithmic}[1]
  \Procedure{eigenvectors}{$\M[\bar]{S},R_n,n$}
  \State $\texttt{myProcID} = (p_1,p_2,\dots,p_N)$
  \State $\texttt{myProcCol} \gets (p_1, \dots, p_{n-1}, \ast, p_{n+1}, \dots, p_N)$
  \State $\M{S} = \textsc{All-Gather}(\M[\bar]{S},\texttt{myProcCol})$ 
  \label{line:eigenvectors:all-gather}
  \State $\Mn{U}{n} = \textsc{Local-Eigenvectors}(\M{S},R_n)$ 
  \label{line:eigenvectors:dsyev}
  \State $\Mn[\bar]{U}{n} = \textsc{Row-Subset}(\Mn{U}{n},P_n,p_n)$
  \Comment Extract $p_n$-th block
  \label{line:eigenvectors:extract}
  \State\Return $\Mn[\bar]{U}{n}$
\EndProcedure
\end{algorithmic}
\end{algorithm}

The cost and memory of the eigenvector computation (\cref{alg:eigenvectors})  is
\begin{align*}
  C_{\textsc{eig}} = &
  \underbrace{\alpha \log P_n + \beta \frac{P_n-1}{P_n} I_n}_{\text{\cref{line:eigenvectors:all-gather}}} + 
  \underbrace{\gamma \frac{10}{3} I_n^3}_{\text{\cref{line:eigenvectors:dsyev}}}, \text{ and}\\
  M_{\textsc{eig}} =  &
  \underbrace{I_n^2/P_n}_{\M[\bar]{S}} +
  \underbrace{I_n^2}_{\M{S}} +
  \underbrace{I_nR_n}_{\Mn{U}{n}} +
  \underbrace{I_nR_n/P_n}_{\Mn[\bar]{U}{n}}.
\end{align*}
The memory for both ``large'' matrices ($\M{S}$ and $\Mn{U}{n}$) is temporary and can overlap the local portions.

\section{Parallel Algorithm Analysis}
\label{sec:algorithms}

\subsection{ST-HOSVD}
\label{sec:st-hosvd}

The initialization of the factor matrices in \cref{alg:hooi} is computed using ST-HOSVD in \cref{alg:sthosvd}. 
The latency cost for ST-HOSVD is dominated by the TTMs, which is given by 
\begin{inlinemath}
\alpha \sum_{n=1}^N  P_n \log P_n.
\end{inlinemath}
The remaining bandwidth and flop contribution of TTM within ST-HOSVD is 
\begin{equation}\label{eq:ttm-sthosvd}\nonumber
\frac{1}{P} \sum_{n=1}^N \lt[
(\beta (P_n-1) + 2 \gamma I_n )
\prod_{k \leq n} R_k \prod _{k > n} I_k 
\rt].
\end{equation}
The total contribution of Gram to the %
cost of ST-HOSVD (aside from latency) is
\begin{multline}\label{eq:gram-sthosvd}\nonumber
  \frac 1P \sum_{n=1}^N \biggl[
  (2 \gamma I_n + 2(P_n-1)\beta) \prod_{k<n} R_k \prod_{k\geq n} I_k + 
   2 \beta (\hat P_n -1 ) I_n^2 
  \biggr] %
\end{multline}
The total contribution of calculating eigenvectors to the cost of ST-HOSVD (aside from latency) is 
\begin{displaymath}
  \sum_{n=1}^N \biggl[ \beta \frac{P_n-1}{P_n} I_n^2 + \gamma \frac{10}{3} I_n^3
  \biggr].
\end{displaymath}

The eigenvector cost is typically negligible compared to TTM and Gram.
Compared to TTM, Gram has a factor of 2 on the bandwidth cost as well as a subtle change in limits of the product notation: 
the $n$th term in the Gram summation is a factor of $I_{n}/R_{n}$ larger than the $n$th term in the TTM summation. 

Note that the dominant expense in TTM and Gram depends on the ordering of modes; each iteration reduces the working data size. 
We can permute the modes in the main iteration to optimize the total cost (both computation and communication). 
We show some examples on the impact of these rearrangements in \cref{sec:vary_mode_order}.

In terms of memory, we can reuse the local variables from iteration to iteration and need not maintain temporaries. 
Therefore, the maximum storage \emph{per processor} for ST-HOSVD is bounded above by
\begin{equation}\label{eq:memory}
  2I/P + \sum_{n=1}^N R_nI_n/P_n + \max_n I_n^2 + \max_n R_nI_n.
\end{equation}

\subsection{HOOI}
\label{sec:hooi}

The iterative improvement of the solution is computed using HOOI in \cref{alg:hooi}.
We derive the cost per outer iteration.
The $N$ inner iterations (computing each factor matrix update) can be done in any order, and this order does not affect the outer iteration cost.
The latency cost of the TTMs dominates the cost of HOOI and each iteration is bounded above by
\begin{inlinemath}
  \alpha N \sum_{n=1}^N P_n \log P_n .
\end{inlinemath}
The bandwidth and computation cost from all $N(N-1)$ TTMs within one outer iteration of HOOI (ignoring the final TTM) is
\begin{displaymath}
  \frac{N-1}{P}  \sum_{n=1}^N (\beta (P_n-1) + 2 \gamma I_n) 
     \prod_{k \leq n} R_{k} \prod_{k>n} I_{k}.
\end{displaymath}
The total contribution of Gram (aside from latency) to the cost of one outer iteration of HOOI is
\begin{equation*}
\frac1P \sum_{n=1}^N (2\beta(P_n-1) + 2\gamma I_n)I_n\hat R_n + 2\beta I_n^2(\hat P_n-1).
\end{equation*}
The total contribution of the eigenvector calculations (aside from latency) to the cost of one outer iteration of HOOI is
\begin{displaymath}
  \sum_{n=1}^N \lt( \beta \frac{P_n-1}{P_n} I_n + \gamma \frac{10}{3} I_n^3 \rt).
\end{displaymath}
This cost is typically dominated by the multiple-TTM computations in \cref{line:hooi:Y}, particularly so by the first TTM performed in each mode.
As in ST-HOSVD, the multiple-TTM computations can be performed in any order, which can greatly affect run time.
In fact, each of the $N$ multiple-TTM computations can be ordered independently.
We do not tune over these possibilities in this work.

The maximum memory is bounded above in the same way as for ST-HOSVD; see \cref{eq:memory}.

\section{Application to Direct Numerical Simulation Data in Combustion
Science}
\label{sec:combustion}

We demonstrate the utility of Tucker compression for two scientific data sets obtained by direct numerical simulation (DNS).
A single simulation today can easily produce 100-1000~GB of data, and much more is expected in the future. 
Some attempts at compression using PCA and other methods have been made; see, e.g., \cite{YaPoCh13}.
Current data sizes are an obstacle for visualization and analysis because the data is difficult to transfer and requires high-end workstations or parallel clusters for computation. 
The goal of compression is to enable easier sharing of data and to facilitate analysis on reconstructed portions of the data. 
For instance, without reconstructing the entire data set, we can extract only the reconstruction of a single species, a few time steps, a coarser grid, a subset of the grid, or any combination of these.
This enables the same data analysis to be performed on laptops.

\subsection{Data Description}
\label{sec:data-description}

The DNS code used to produce this data is called S3D, a massively parallel compressible reacting flow solver developed at Sandia National Laboratories \cite{ChChSuDe09}.
We work with the following multiway data sets:
\begin{itemize}%
\item \textbf{HCCI:} 
  This 4-way data tensor of size $672 \times 672 \times 33 \times 627$ comes from a $672 \times 672$ spatial grid with 33 variables over 627 time steps. Each time step requires 111~MB of storage, and the entire dataset is 70~GB.
  This data comes from the simulation of an autoignitive premixture of air and ethanol in Homogeneous Charge Compression Ignition (HCCI) mode~\cite{BhChLu14}. 
The first two dimensions correspond to the 2D spatial grid, the third to the species, and the last to time.
\item \textbf{TJLR:} This 5-way data tensor of size $460 \times 700 \times 360 \times 35 \times 16$ comes from a $460 \times 700 \times 360$ spatial grid with 35 variables (30 species plus five derivative values) over 16 time steps. Each time step requires 32~GB storage, so the entire dataset is 520~GB.
  It comes from a temporally-evolving planar slot jet flame with DME (dimethyl ether) as the fuel \cite{BhLuShSu15}. 
  This data is low resolution due to  significant downsampling and 
so  is less amenable to compression than the HCCI dataset.
\item \textbf{SP:} This 5-way data tensor is of size $500 \times 500 \times 500 \times 11 \times 50$ and corresponds to a cubic $500\times500\times500$ spatial grid for 11 variables over 50 time steps.  Each time step requires 11~GB, so the entire dataset is 550~GB.  
The SP dataset is from the simulation of a three-dimensional statistically steady planar turbulent premixed flame of methane-air combustion. The simulation considers reduced chemical kinetics of 6 species, which results in a total of 11 simulation variables \cite{KZCS16}.  
Unlike the HCCI and TJLR datasets, the SP dataset represents a statistically steady, rather than a temporally evolving, turbulent flame. This means the SP data set is more amenable to compression.
\end{itemize}
Each data set is centered and scaled \emph{for each variable/species}. 
We compute the mean and standard deviation for each species slice, and then we transform the data by subtracting the mean and dividing the result by the standard deviation (unless it is less than $10^{-10}$, in which case the division is not performed).
This normalizes the data so that we can roughly assume that each entry comes from a standard normal distribution.

\newcommand{\FigModeErrors}[3][north east]{
  \begin{tikzpicture}[scale=.7]
    \renewcommand{\datafile}{#2}
    \pgfmathparse{1e-3/sqrt(#3)}
    \edef\err{\pgfmathresult}
    \begin{semilogyaxis}[
      legend pos = #1,
      ylabel={Mode-wise Norm. RMS Error}, 
      xlabel={Rank},
      unbounded coords=jump,
      y tick label style={rotate=90},
      ylabel near ticks,
      ]	

      \ifdimless{#3 pt}{5 pt}
      {\legend{Spatial 1, Spatial 2, Species, Time}}
      {\legend{Spatial 1, Spatial 2, Spatial 3, Species, Time}}

      \foreach \foo in {1,...,#3}
      \addplot table[x=Rank, y=M\foo] {\datafile};

      \draw[dashed] (axis cs:\pgfkeysvalueof{/pgfplots/xmin},\err) -- (axis cs:\pgfkeysvalueof{/pgfplots/xmax},\err);
      \node [below left] at (axis cs:\pgfkeysvalueof{/pgfplots/xmax},\err) {$\epsilon/\sqrt{N}$};
      
    \end{semilogyaxis}
  \end{tikzpicture}
}
%
%
%
%
\begin{figure*}[tb]
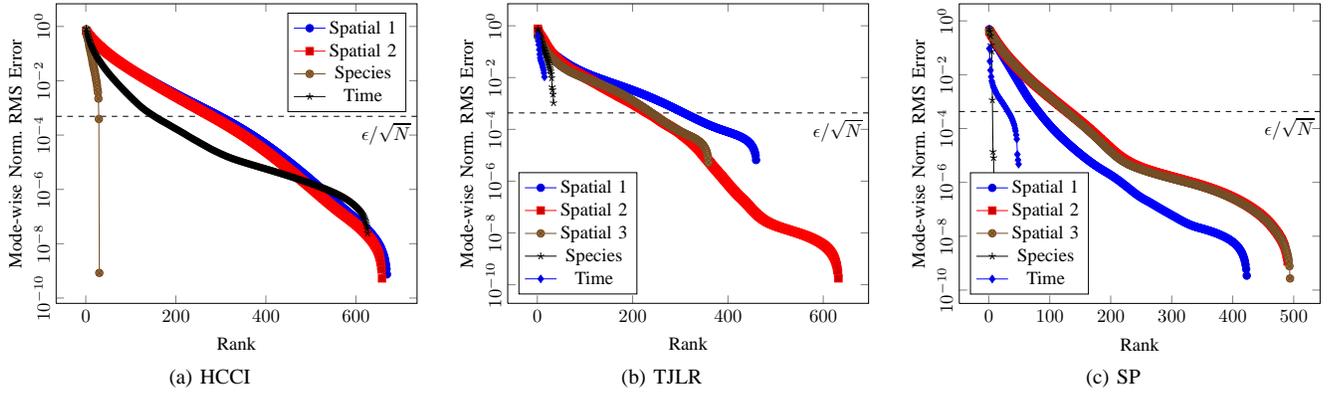

\subfloat[HCCI]{
\FigModeErrors{mode_errors_hcci_627.dat}{4}
} 
\subfloat[TJLR]{
\FigModeErrors[south west]{mode_errors_tjlr_16.dat}{5}
}
\subfloat[SP]{
\FigModeErrors[south west]{mode_errors_sp_50.dat}{5}
}
\caption{Mode-wise contributions to the error bound for combustion data sets.}
\label{fig:mode_errors}
\end{figure*}

\begin{table*}
  \centering\footnotesize
  \caption{Compression and errors for maximum normalized RMS error threshold of 1e-3.}
  \label{tab:errors}
  \begin{tabular}{| c | c | c | c | c | c | c |}
    \hline
    \multirow{3}{*}{Dataset} & \multirow{3}{*}{Reduced Dimensions} & \multicolumn{2}{c|}{ST-HOSVD} & \multicolumn{2}{c|}{HOOI} & \multirow{2}{*}{Compression} \\
    & & Norm. & Max. Abs. & Norm. & Max. Abs. & \multirow{2}{*}{Ratio} \\
    & & RMS & Elem. Err. & RMS & Elem. Err. & \\
    \hline
    HCCI & ($297$, $279$, $29$, $153$) & 9.259e-4 & 1.510e-1 & 9.254e-4 & 1.509e-1 & 25  \\
    TJLR  & ($306$, $232$, $239$, $35$, $16$) & 7.617e-4 & 1.567 & 7.617-4 & 1.568 & 7 \\
    SP  & ($81$, $129$, $127$, $7$, $32$) & 8.663e-4 & 1.017e-1 & 8.662e-4 & 1.017e-1 & 231 \\
    \hline
  \end{tabular}
\end{table*}

\subsection{Compression Rates}
\label{sec:compression-rates}

Let $\T{X}$ be an $N$-way data tensor of size $I_1 \times \cdots \times I_N$. 
Let $\lambda_i^{(n)}$ denote the $i$th eigenvalue of the Gram matrix $\Mz{X}{n}\Mz{X}{n}^{\Tra}$ for $n=1,\dots,N$, corresponding to the square of the $i$th singular value of $\Mz{X}{n}$.
Assume that the eigenvalues are in decreasing order, i.e., $\lambda_1^{(n)} \geq \lambda_2^{(n)} \geq \cdots \geq \lambda_{I_n}^{(n)}$.
Let $\T[\tilde]{X}$ be the reconstruction per \cref{eq:reconstruct} using the T-HOSVD with reduced size $R_1 \times \cdots \times R_N$.
It is well known \cite{DeDeVa00,VaVaMe12} that 
selecting $R_n$ such that
\begin{displaymath}
  \sum_{i=R_n+1}^{I_n} \lambda_i^{(n)} \leq \epsilon^2 \|\T{X}\|^2 / N 
\end{displaymath}
for $n=1,\dots,N$ ensures that the T-HOSVD satisfies
\begin{equation}
\label{eq:error_bound}
  \| \T{X} - \T[\tilde]{X} \|^2  \leq 
  \sum_{n=1}^N \left( \sum_{i=R_n+1}^{I_n} \lambda_i^{(n)} \right) \leq \epsilon \| \T{X}\|.
\end{equation}
The normalized RMS error for HOOI initialized by ST-HOSVD is bounded above by the T-HOSVD error \cite{DeDeVa00,VaVaMe12}.
In practice, computing $\{R_n\}$ given $\epsilon$ can be done within the ST-HOSVD as specified in \cref{alg:sthosvd} (the T-HOSVD need not be computed).

\Cref{fig:mode_errors} shows the normalized mode-wise RMS, i.e., 
\begin{displaymath}
  \lt(\sum_{i=R_n+1}^{I_n} \lambda_i^{(n)}\rt)^{1/2} / \|\T{X}\|
\end{displaymath}
for each mode and value of $R_n$ of HCCI-628. 
The rate of drop-off in the errors determines the compressibility of the data, as the intersections of the mode-wise error curves with a desired threshold (shown by the dotted line) give upper bounds on the reduced dimensions for $\epsilon=10^{-3}$.

We show compression rates for all data sets in \cref{fig:compression}, where
the compression ratio is
\begin{displaymath}
  C = {\prod_{k=1}^N I_n} \bigg/ \lt( {\prod_{k=1}^N R_n + \sum_{k=1}^N I_n R_n} \rt).
\end{displaymath}
The  TJLR data set is the least compressible data set with $C$ ranging from 2 at  $\epsilon=10^{-6}$ to 37 for $\epsilon=10^{-2}$. 
The SP data set is much more compressible, with $C$ ranging from 5 to 5600 for the same error range.

\begin{figure}[htb]
\centering
\begin{tikzpicture}[scale=0.8]
\renewcommand{\datafile}{compression.dat}

\begin{loglogaxis}[
	legend pos = north west,
	ylabel={Compression Ratio}, 
	xlabel={Max.\@ Normalized RMS Error},
	y tick label style={rotate=90},
        ylabel near ticks,
]
	\addplot table[x=Error, y=hcci_627] {\datafile}; \addlegendentry{HCCI}
	\addplot table[x=Error, y=tjlr_16] {\datafile}; \addlegendentry{TJLR}
	\addplot table[x=Error, y=sp_50] {\datafile}; \addlegendentry{SP}
\end{loglogaxis}
\end{tikzpicture}
%
%
%
%

\caption{Approximation error versus compression for different data sets.}
\label{fig:compression}
\end{figure}
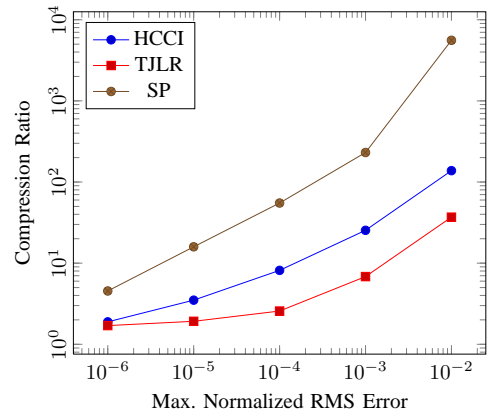

\begin{figure*}
 \centering
\subfloat[Varying processor grid for tensor of size $384{\times}384{\times}384{\times}384$ with reduced size of $96{\times}96{\times}96{\times}96$.]
{
\hspace{0.25in}
 \label{fig:vary_proc_grid}
\begin{tikzpicture}[scale=.8]
\renewcommand{\datafile}{vary_proc_grid.dat}
\pgfplotstablesort[sort key={Total},sort cmp={float <}]{\sorted}{\datafile}
\pgfplotstablegetelem{0}{Total}\of{\sorted}
\let\mintot=\pgfplotsretval
\def\xlistmacro{}
\def\xliststring{}
\pgfplotstableforeachcolumnelement{ProcGrid}\of\sorted\as\entry{%
\xifinlist{\entry}{\xlistmacro}{}{
        \listxadd{\xlistmacro}{\entry}
        \edef\xliststring{\xliststring\entry,}
    }
}
\begin{axis}[
  	height=3in,
	ybar stacked,
	legend pos = north west,
	reverse legend,
	ymin=0,
	xtick=data,
	symbolic x coords/.expand once={\xliststring},
	x tick label style={rotate=60,anchor=east},
]

	\pgfplotsset{cycle list={red, fill=red \\ green, fill=green \\ blue, fill=blue \\}}

	\legend{Gram, Evecs, TTM}

	\addplot table[x=ProcGrid, y expr=\thisrow{Gram-0}/\mintot] {\datafile};       
	\addplot table[x=ProcGrid, y expr=\thisrow{Evec-0}/\mintot] {\datafile};
	\addplot table[x=ProcGrid, y expr=\thisrow{TTM-0}/\mintot] {\datafile};
	\addplot table[x=ProcGrid, y expr=\thisrow{Gram-1}/\mintot] {\datafile};
	\addplot table[x=ProcGrid, y expr=\thisrow{Evec-1}/\mintot] {\datafile};
	\addplot table[x=ProcGrid, y expr=\thisrow{TTM-1}/\mintot] {\datafile};
	\addplot table[x=ProcGrid, y expr=\thisrow{Gram-2}/\mintot] {\datafile};
	\addplot table[x=ProcGrid, y expr=\thisrow{Evec-2}/\mintot] {\datafile};
	\addplot table[x=ProcGrid, y expr=\thisrow{TTM-2}/\mintot] {\datafile};
	\addplot table[x=ProcGrid, y expr=\thisrow{Gram-3}/\mintot] {\datafile};
	\addplot table[x=ProcGrid, y expr=\thisrow{Evec-3}/\mintot] {\datafile};
	\addplot %
        table[x=ProcGrid, y expr=\thisrow{TTM-3}/\mintot] {\datafile};

\end{axis}
\end{tikzpicture}

%
%
%
%
\hspace{0.25in}
} 
\hspace{0.25in}
\subfloat[Varying mode order for tensor of size $25{\times}250{\times}250{\times}250$ with reduced size $10{\times}10{\times}100{\times}100$.]
{
\hspace{0.25in}
  \label{fig:vary_mode_order}
\begin{tikzpicture}[scale=0.8]

\renewcommand{\datafile}{vary_mode_order.dat}

\begin{axis}[
  	height=3in,
	ybar stacked,
	legend pos = north west,
	reverse legend,
	ymin=0,
	xtick=data,
	symbolic x coords={1023,1203,1230,0123,0213,0231,2013,2031,2103,2130,2301,2310},
	xticklabels={~~~~~1234,1324,1342,2134,2314,2341,3124,3142,3214,3241,3412,3421},
	x tick label style={rotate=60,anchor=east},
]
	\pgfplotstablesort[sort key={Total},sort cmp={float <}]{\sorted}{\datafile}
    	\pgfplotstablegetelem{0}{Total}\of{\sorted}
    	\let\mintot=\pgfplotsretval

	\pgfplotsset{cycle list={red, fill=red \\ green, fill=green \\ blue, fill=blue \\}}

	\legend{Gram, Evecs, TTM}

	\addplot table[x=Perm, y expr=\thisrow{Gram-0}/\mintot] {\datafile};
	\addplot table[x=Perm, y expr=\thisrow{Evec-0}/\mintot] {\datafile};
	\addplot table[x=Perm, y expr=\thisrow{TTM-0}/\mintot] {\datafile};
	\addplot table[x=Perm, y expr=\thisrow{Gram-1}/\mintot] {\datafile};
	\addplot table[x=Perm, y expr=\thisrow{Evec-1}/\mintot] {\datafile};
	\addplot table[x=Perm, y expr=\thisrow{TTM-1}/\mintot] {\datafile};
	\addplot table[x=Perm, y expr=\thisrow{Gram-2}/\mintot] {\datafile};
	\addplot table[x=Perm, y expr=\thisrow{Evec-2}/\mintot] {\datafile};
	\addplot table[x=Perm, y expr=\thisrow{TTM-2}/\mintot] {\datafile};
	\addplot table[x=Perm, y expr=\thisrow{Gram-3}/\mintot] {\datafile};
	\addplot table[x=Perm, y expr=\thisrow{Evec-3}/\mintot] {\datafile};
	\addplot table[x=Perm, y expr=\thisrow{TTM-3}/\mintot] {\datafile};

\end{axis}
\end{tikzpicture}

%
%
%
%
\hspace{0.25in}
}
\caption{Relative run-time comparisons for different tuning options in ST-HOSVD.}
\end{figure*}
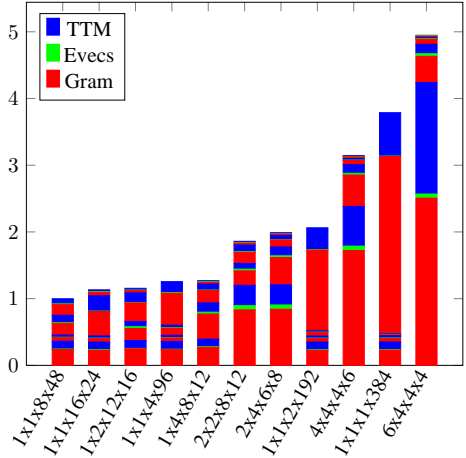
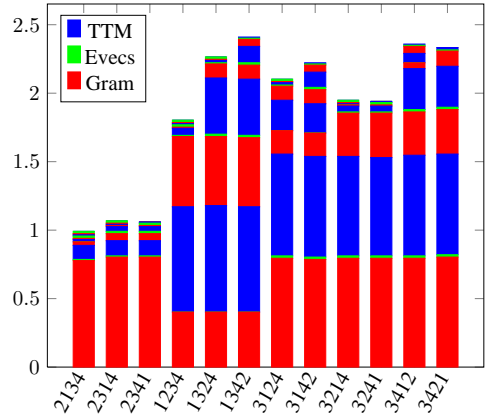

\subsection{Reconstruction Error}
\label{sec:reconstruction-error}

\Cref{tab:errors} presents detailed compression results using $\epsilon = 10^{-3}$ to determine the compression ratios. We include the maximum absolute element error of the centered and scaled data.
The HOOI iterations make little improvements on the ST-HOSVD initialization, so simply performing ST-HOSVD (with no HOOI iterations) is likely sufficient for this particular application area.
  
\section{Performance Results}
\label{sec:performance}

\subsection{Experimental Platform}

We run all experiments on Edison, a Cray XC30 supercomputer located at NERSC consisting of 5,576 dual-socket 12-core Intel ``Ivy Bridge'' (2.4 GHz) compute nodes. 
The peak flop rate of each core is 19.2 GFLOPS. %
Each node has 64 GB of memory. %
The nodes are connected by a Cray ``Aries'' interconnect with a dragonfly topology.
We use Cray compilers and LibSci for BLAS and LAPACK subroutines.

\subsection{Parameter Choice: Processor Grid Configuration}
\label{sec:vary_proc_grid}

Our first microbenchmark demonstrates the effect of the processor grid on performance of the ST-HOSVD algorithm.
\Cref{fig:vary_proc_grid} presents relative running times for a fixed problem size and number of processors, varying only the processor grid.
We break down the running time across the three subroutines---Gram, Evecs, and TTM---which are each performed once per mode.
Thus, each bar has four blocks of each color, ordered from bottom to top so that the first subroutine computation corresponds to the bottommost block of that color.
The dimensions are all equal to ensure the ordering of the modes does not impact the performance.

As shown in \cref{sec:st-hosvd}, the processor grid does not change the number of flops for each step in ST-HOSVD, but
it does affect the performance of sequential linear algebra kernels due to its impact on the dimensions and layouts of the local matrices involved in the computations.

For most processor grids, the initial iteration consumes at least half of the overall running time. 
Gram dominates this initial iteration.
As shown in the analysis of \cref{sec:st-hosvd}, the first Gram is more expensive than the first TTM by a factor of at least $I_1/R_1=4$ (in terms of both computation and communication).
The cost of the Eigenvectors computation is negligible.
The best processor grids have $P_1=1$, so that communication is minimized in the first iteration. In this case, the first Gram calculation only calls the all-reduce, which is of size $I_1^2$, and
the first TTM involves no communication at all.
We do not show results for processor grids with $P_1>6$, as they are more than 5 times the optimal running time.

The best processor grid for ST-HOSVD is not necessarily optimal for HOOI.
While our limited tuning suggests that good choices for ST-HOSVD tend to be reasonable for HOOI and vice versa, an optimal overall processor grid choice for \cref{alg:hooi} depends on the number of iterations of HOOI.

\begin{figure*}
\centering
\subfloat[Strong scaling for $200{\times}200{\times}200{\times}200$ tensor with reduced size $20{\times}20{\times}20{\times}20$, using $2^k$ nodes for $0\leq k\leq 9$.]{
\centering
\begin{tikzpicture}[scale=.8]
\renewcommand{\datafile}{strong_scaling.dat}
\begin{axis}[
	legend pos = north east,
	reverse legend,
	ylabel={Time (seconds)}, 
	xlabel={Number of Nodes},
	xtick=data,
	xticklabels={1,2,4,8,16,32,64,128,256,512},
	xmode=log,
	log basis x ={2},
	ymode=log,
	log basis y={2},
	y tick label style={rotate=90},
        ylabel near ticks,
]
	\legend{ST-HOSVD,HOOI}

	\addplot table[x expr=\thisrow{Cores}/24, y=HOSVD-s] {\datafile};
	\addplot table[x expr=\thisrow{Cores}/24, y=HOOI-s] {\datafile};

\end{axis}
\end{tikzpicture}
%
%
%
%
\label{fig:strong_scaling}
} \hspace{1in}
\subfloat[Weak scaling for $200k{\times}200k{\times}200k{\times}200k$ tensor with reduced size $20k{\times}20k{\times}20k{\times}20k$, using $k^4$ nodes for $1\leq k\leq 6$.]{
\centering
\begin{tikzpicture}[scale=.8]
\renewcommand{\datafile}{weak_scaling.dat}

\begin{axis}[
	legend pos = north east,
	ylabel={GFLOPS Per Core}, 
	xlabel={Number of Nodes},
	xtick=data,
	xticklabels={1,,81,256,625,1296},
	ymin=0,
	ymax=19.2,
	ytick={5,10,15,19.2},
	y tick label style={rotate=90},
        ylabel near ticks,
        reverse legend
]
	\legend{ST-HOSVD,HOOI}

	\addplot table[x expr=\thisrow{Cores}/24, y=HOSVD-G] {\datafile};
	\addplot table[x expr=\thisrow{Cores}/24, y=HOOI-G] {\datafile};

\end{axis}
\end{tikzpicture}
\label{fig:weak_scaling}
}
\caption{Scaling performance of ST-HOSVD and one iteration of HOOI.}
\end{figure*}
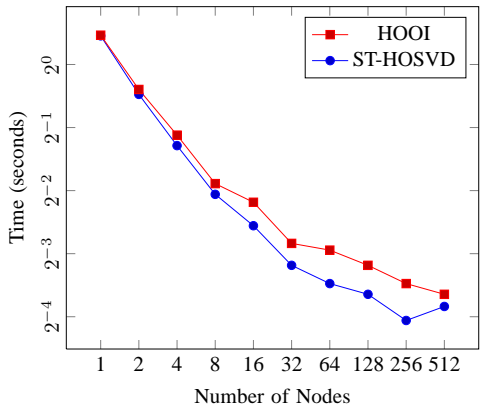
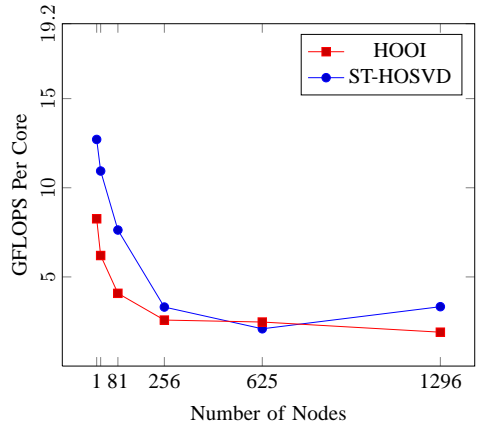

\subsection{Parameter Choice: Mode Ordering}
\label{sec:vary_mode_order}

Our second microbenchmark demonstrates the effect of the mode ordering on performance of the ST-HOSVD algorithm.
In \cref{fig:vary_mode_order}, we show relative run times for a fixed problem size and processor grid, varying only the order of modes in \cref{line:hosvd:loop} of \cref{alg:sthosvd}.
As in \cref{fig:vary_proc_grid}, we break down the run times across the three subroutines so that each bar comprises a block of each subroutine color for each mode.

For this problem, we use synthetic data. The full tensor dimensions are $25\times250\times250\times250$ which is formed from a Tucker decomposition with core dimensions $10\times10\times100\times100$. We use a processor grid of dimensions $2\times2\times2\times2$.
We vary both the tensor dimensions and the ratios of tensor to core dimensions to accentuate the effects of mode ordering.
We use only 16 of 24 cores on one node in order to obtain a uniform processor grid and eliminate variability based on processor grid choice.

Most of the overall performance is determined by the choice of the first dimension, which is the smallest by a factor of 10.
Starting with the first dimension means a cheaper first Gram but less reduction in computation and communication in subsequent iterations.
Starting with the second dimension, which has the largest compression ratio, yields the greatest savings in subsequent iterations but incurs more overhead in the first iteration.
Nevertheless,
the optimal ordering starts with the second dimension.
While the first Gram computation is more expensive than starting with the first dimension, the reduction in time for subsequent operations more than compensates.

We note that the authors of ST-HOSVD propose a heuristic ordering (in the case of sequential computation) that greedily chooses the mode that minimizes the number of flops in the current iteration \cite{VaVaMe12}.
While that heuristic is not optimal in this case, we see no simple alternative scheme that is always optimal.
Another reasonable heuristic is to greedily choose the mode that maximizes the compression ratio $I_n/R_n$.

\subsection{Strong Scaling}
\label{sec:strong_scaling}

To test the parallel scaling of our algorithm, we fix a particular problem and increase the number of processors we use to compute the Tucker decomposition.
\Cref{fig:strong_scaling} reports the running time of ST-HOSVD and one iteration of HOOI across $24\cdot 2^k$ processors, for $0\leq k\leq 9$.

For the experiment, we set the tensor dimensions to be $200\times200\times200\times200$ and compress it to a core tensor of dimension $20\times20\times20\times20$.
For each $k$, we tune the processor grid over three or four possibilities (chosen heuristically based on the results of \cref{sec:vary_proc_grid}) for each algorithm and report the minimum running time.

On one node, the data set requires about 1/5th of the available memory.
ST-HOSVD and one step of HOOI runs in about 3 seconds, achieving 67\% and 45\% of peak performance, respectively.
On 512 nodes, the same computation requires 0.15 seconds, achieving aggregate performance of 6.0 and 3.8 TFLOPS, respectively.
At this scale, the local tensor data is about 1 MB per core.
We note that we continue to decrease running time up to 256 nodes (6144 cores), though we are far from peak performance, as the interprocessor communication and small matrix dimensions within local computation kernels both degrade performance.

\subsection{Weak Scaling}
\label{sec:weak_scaling}

\Cref{fig:weak_scaling} shows a weak scaling experiment, reporting performance per core for both ST-HOSVD and one iteration of HOOI.
In this experiment, we fix the amount of data per processor and increase the number of processors and tensor dimensions simultaneously.

We use $24\cdot k^4$ processors and set the tensor dimensions to be $(200k)^4$ with cores of dimension $(20k)^4$, for $1\leq k\leq 6$.
Note that $k=1$ corresponds to the dimensions we use in the strong scaling experiment in \cref{sec:strong_scaling}.
The plot reports the best performance for each algorithm over three different processor grids: $1\times 1\times 4k^2\times 6k^2$, $k\times k\times 4k\times 6k$, and $k\times 2k\times 3k\times 4k$.
The size of the data sets ranges from about 12~GB for $k=1$ up to 15~TB for $k=6$.
On one node ($k=1$),
ST-HOSVD and HOOI achieve 66\% and 43\% of peak performance, respectively.
For 1296 nodes ($k=6$), the algorithms achieves 17\% and 12\% of peak, respectively, or up to 104 TFLOPS in aggregate.
The time required to process the 15~TB data set (performing ST-HOSVD and HOOI on data in memory) is 70 seconds.

The main reason that we see degradation in performance as we scale up to high processor counts is that it becomes harder to navigate the tradeoffs among optimization of the processor grid for the computations in different modes.

\section{Conclusion}
\label{sec:conclusion}

Our parallel Tucker decomposition enables compression of massive data sets that do not fit into memory on a single machine. We show that these data sets can be processed in reasonable time and yield excellent compression rates.
The implementation performs well (near peak performance) at low core counts, and it scales (offers reduced run times) up to high core counts.
In order to achieve even better parallel scaling, we see several avenues for performance improvement
such as using multi-threaded {BLAS} for all local computations or adapting
recent work in optimizing multi-threaded TTM computations \cite{LBPSV15}.
Additionally, we can overlap communication and computation and fully exploit the symmetry in the Gram computation. %
For achieving approximation errors near the square root of machine precision (or smaller), we need to consider a numerical improvement to our algorithm, directly computing the singular values.
Improving the numerical stability of our algorithm will not drastically hurt our overall performance;
because $\Mz{Y}{n}^{\Tra}$ is typically very tall and skinny, we can compute the SVD using a QR decomposition as a preprocessing step at roughly twice the cost of our current approach.
\section*{Acknowledgment}
\small
We thank Hemanth Kolla and Ankit Bhagatwala for providing the
combustion application data.
This research used resources of the National Energy Research Scientific Computing Center and the Oak Ridge Leadership Computing
Facility, DOE Office of Science User Facilities supported by the Office of Science of the U.S. Department of Energy under Contract Nos. DE-AC02-05CH11231 and DE-AC05-00OR22725, respectively.
This material is based upon work supported by the Sandia Truman
Postdoctoral Fellowship 
and the U.S. Department of
Energy, Office of Science, Office of Advanced Scientific Computing
Research, Applied Mathematics program.
Sandia National Laboratories is a multi-program laboratory managed and
operated by Sandia Corporation, a wholly owned subsidiary of Lockheed
Martin Corporation, for the U.S. Department of Energy's National
Nuclear Security Administration under contract DE--AC04--94AL85000.
\bibliographystyle{siamplain}

%

%
%
%
\end{document}